\documentclass{ws-procs9x6}
\usepackage{amssymb}
\usepackage{amsfonts}
\usepackage{amscd}
\usepackage{latexsym}
\usepackage{enumerate}
\usepackage{amsmath}
\usepackage{xy} \xyoption{all}

\numberwithin{equation}{section}

\newcounter{cs}
\stepcounter{cs}
\newcounter{ds}
\stepcounter{ds}

\newcommand{\casos}{\begin{itemize}}
\newcommand{\fcasos}{\end{itemize}\setcounter{cs}{1}}

\newcommand{\ol}{\overline}

\newcommand{\Id}{{\mbox{Id}}}

\newcommand{\PP}{\mathbb{P}}
\newcommand{\Pfree}{\PP_{\mathrm{free}}}
\newcommand{\Preg}{\PP_{\mathrm{reg}}}

\DeclareMathOperator{\rL}{L}
\newcommand{\rLfree}{\rL_{\mathrm{free}}}
\newcommand{\rLreg}{\rL_{\mathrm{reg}}}

\newcommand{\dnw}{\mathbin{\downarrow}}

\newcommand{\R}{\mathbb{R}}

     \newcommand{\mon}[1]{\mathcal{V}(#1)}              
     
\begin{document}

\title{\large\sf THE REALIZATION PROBLEM FOR VON NEUMANN\\REGULAR RINGS}

\author{PERE ARA}

\address{Departament de Matem\`atiques, Universitat Aut\`onoma de
Barcelona,\\ 08193, Bellaterra (Barcelona),
Spain\\E-mail: para@mat.uab.cat}

\begin{abstract}
We survey recent progress on the realization problem for von Neumann regular rings,
which asks whether every countable conical refinement monoid can be realized as the monoid of
isoclasses of finitely generated projective right $R$-modules over a von Neumann regular ring $R$.
\end{abstract}

\keywords{von Neumann
regular ring; Leavitt path algebra; refinement monoid.}

\bodymatter


\section*{}

\noindent This survey consists of four sections. Section 1 introduces the realization problem
for von Neumann regular rings, and points out its relationship with the
separativity problem of [\refcite{AGOP}]. Section 2 surveys positive realization results
for countable conical refinement monoids, including the recent constructions in [\refcite{AB}]
and [\refcite{Aposet}]. We analyze in Section 3 the relationship with the realization
problem of algebraic distributive lattices as lattices of two-sided ideals over von Neumann
regular rings. Finally we observe in Section 4 that there are countable conical
monoids which can be realized by a von Neumann regular $K$-algebra for some countable field $K$,
but they cannot be realized by a von Neumann regular $F$-algebra for any uncountable field $F$.

\section{The problem}
\label{sect:prob}

All rings considered in this paper will be associative, and all the monoids
will be commutative.

For a unital ring $R$, let $\mon{R}$ denote the monoid
of isomorphism classes of finitely generated projective right
$R$-modules, where the operation is defined by $$[P]+[Q]=[P\oplus
Q].$$
This monoid describes faithfully the decomposition structure of
finitely generated projective modules. The monoid $\mon{R}$ is
always a {\it conical monoid}, that is, whenever $x+y=0$, we have
$x=y=0$. Recall that an {\it order-unit} in a monoid $M$ is an element
$u$ in $M$ such that for every $x\in M$ there is $y\in M$ and $n\ge 1$
such that $x+y=nu$. Observe that $[R]$ is a canonical order-unit in $\mon{R}$.
By results of Bergman \cite[Theorems 6.2 and 6.4]{B} and
Bergman and Dicks \cite[page 315]{BD}, any conical monoid with an
order-unit appears as $\mathcal V (R)$ for some unital hereditary
ring $R$.

A monoid $M$ is said to be a {\it refinement monoid} in
case any equality $x_1+x_2=y_1+y_2$ admits a refinement, that is,
there are $z_{ij}$, $1\le i,j \le 2$ such that $x_i=z_{i1}+z_{i2}$
and $y_j=z_{1j}+z_{2j}$ for all $i,j$, see e.g. [\refcite{AMP}]. If $R$
is a von Neumann regular ring, then the monoid $\mon{R}$ is a
refinement monoid by \cite[Theorem 2.8]{vnrr}.

\medskip

The following is still an open problem:

\medskip

\noindent {\bf R1. Realization Problem for von Neumann Regular Rings} Is every
 countable conical refinement monoid realizable by a von Neumann
regular ring?

\medskip

A related problem was posed by K.R. Goodearl in [\refcite{directsum}]:

\medskip

\noindent {\bf FUNDAMENTAL OPEN PROBLEM} Which monoids arise as
$\mon{R}$'s for a von Neumann regular ring $R$?

\medskip

It was shown by Wehrung in [\refcite{Wehisrael}] that there are conical
refinement monoids of size $\aleph _2$ which cannot be realized.
If the size of the monoid is $\aleph _1$ the question is open.
Wehrung's approach is related to Dilworth's Congruence Lattice Problem
(CLP), see Section \ref{sect:CLP}. A solution to the latter problem has recently
appeared in [\refcite{WAdv}].

Problem R1 is related to the separativity problem. A class $\mathcal
C$ of modules is called {\em separative} if for all $A,B\in \mathcal
C$ we have
$$A\oplus A\cong A\oplus B\cong B\oplus B \implies A\cong B.$$
A ring $R$ is {\it separative} if the class $FP(R)$ of all finitely generated
projective right $R$-modules is a separative class.
Separativity is an old concept in semigroup theory, see [\refcite{CP}]. A
commutative semigroup $S$ is called {\em separative} if for all $a,b\in S$ we
have $a+a=a+b=b+b \implies a=b$. An alternative characterization is that
a commutative semigroup is separative if and only if it can be embedded in a product
of monoids of the form $G\sqcup \{\infty \}$, where $G$ is an abelian group.
Clearly a ring $R$ is separative if
and only if $\mon{R}$ is a separative semigroup. Separativity
provides a key to a number of outstanding cancellation problems for
finitely generated projective modules over exchange rings, see
[\refcite{AGOP}].

Outside the class of exchange rings, separativity can easily fail.
In fact it is easy to see that a commutative ring $R$ is separative
if and only if $\mon{R}$ is cancellative. Among exchange rings,
however, separativity seems to be the norm. It is not known whether
there are non-separative exchange rings. This is one of the major
open problems in this area. See [\refcite{sexch}] for some classes of
exchange rings which are known to be separative. We single out the
problem for von Neumann regular rings. (Recall that every von
Neumann regular ring is an exchange ring.)

\medskip

\noindent {\bf SP. Is every von Neumann regular ring separative?}

\medskip

We have $(R1 \text{ has positive answer }) \implies (SP \text{ has a
negative answer })$. To explain why we have to recall results of
Bergman and Wehrung concerning {\em existence} of countable
non-separative conical refinement monoids.

Recall that every monoid $M$ is endowed with a natural pre-order,
the so-called {\it algebraic pre-order}, by $x\le y$ iff there is
$z\in M$ such that $y=x+z$. This is the only order on monoids that
we will consider in this paper. A monoid homomorphism $f\colon M\to M'$ is an
{\it order-embedding} in case $f$ is one-to-one and, for $x,y\in M$, we have
$x\le y$ if and only if $f(x)\le f(y)$.

\begin{proposition}\label{prop:embedding}{\rm(cf. [\refcite{Wehsemforum}])} Let $M$ be a
countable conical monoid. Then there is an
order-embedding of $M$ into a countable conical refinement monoid.
\end{proposition}

Let us apply the above Proposition to the conical monoid
$M$ generated by $a$ with the only relation $2a=3a$. Then
$$a+a=a+(2a)=(2a)+(2a)$$
but $a\ne 2a$ in $M$. By Proposition \ref{prop:embedding} there
exists an order-embedding $M\to M'$, where $M'$ is a countable
conical refinement monoid,  and $M'$ cannot be separative.

Thus if R1 is true we can represent $M'$ as $\mon{R}$ for some von
Neumann regular ring and $R$ will be non-separative.

\section{Known cases}
\label{sect:known}

It turns out that only a few cases of R1 are known. In this section I will
describe the positive realization results of which I am aware.

The first realization result is by now a classical one. Recall that
a monoid $M$ is  said to be {\it unperforated} if, for $x,y\in M$
and $n\ge 1$, the relation $nx\le ny$ implies that $x\le y$. A {\em
dimension monoid} is a cancellative, refinement, unperforated conical monoid.
These are the positive cones of the dimension groups \cite[Chapter
15]{vnrr}. Recall that, by definition, an {\it ultramatricial} $K$-algebra $R$ is
a direct limit of a sequence of finite direct products of matrix
algebras over $K$. Clearly every ultramatricial algebra
is von Neumann regular.

\begin{theorem} {\rm ([\refcite{Elliott}], \cite[Theorem 3.17]{poag}, \cite[Theorem 15.24(b)]{vnrr} )}
\label{thm:dim}
If $M$ is a countable dimension monoid and $K$ is any
field, then there exists an ultramatricial $K$-algebra $R$
such that $\mon{R}\cong M$.
\end{theorem}

A $K$-algebra is said to be {\it locally matricial} in case it is a direct limit of a directed system
of finite direct products of matrix algebras over $K$, see \cite[Section 1]{tensorp}.
It was proved in \cite[Theorem 1.5]{tensorp} that if $M$ is a dimension monoid
of size $\le \aleph _1$, then
it can be realized as $\mon{R}$ for a locally matricial $K$-algebra $R$.
Wehrung constructed in [\refcite{Wehisrael}] dimension monoids of size $\aleph _2$
which cannot be realized by regular rings. Indeed the monoids constructed in
[\refcite{Wehisrael}] are the positive cones of dimension groups which are vector spaces
over $\mathbb Q$. A refinement of the method used in [\refcite{Wehisrael}] gave a
dimension monoid counterexample of size $\aleph _2$ with an order unit of index two [\refcite{WJalg}],
thus answering a question posed by Goodearl in [\refcite{bounded}].

\smallskip

Another realization result was obtained by Goodearl, Pardo and the
author in [\refcite{AGP}].

\begin{theorem}\label{thm:purinf}
\cite[Theorem 8.4]{AGP}. Let $G$ be a countable abelian group and
$K$ any field. Then there is a purely infinite simple regular
$K$-algebra $R$ such that $K_0(R)\cong G$.
\end{theorem}

Recall that a simple ring $R$ is {\it purely infinite} in case it is
not a division ring and, for every nonzero element $a\in R$ there
are $x,y\in R$ such that $xay=1$ (see \cite[Section 1]{AGP},
especially Theorem 1.6). Since $\mon{R}=K_0(R)\sqcup \{0\}$ for a
purely infinite simple regular ring \cite[Corollary 2.2]{AGP}, we
get that all monoids of the form $G\sqcup \{0\}$, where $G$ is a
countable abelian group, can be realized.

As Fred Wehrung has kindly pointed out to me, another class of
conical refinement monoids which can be realized by von Neumann regular rings is the
class of {\it continuous dimension scales}, see \cite[Chapter 3]{GWMem}. All the monoids in this class
satisfy the property that every bounded subset has a supremum, as well as some additional axioms,
see \cite[Definition 3-1.1]{GWMem}.
Indeed, if $M$ is a commutative monoid, then $M\cong \mon{R}$  for some
regular, right self-injective ring $R$ if and only if $M$ is a continuous
dimension scale with order-unit \cite[Corollary 5-3.15]{GWMem}.
These monoids have unrestricted cardinality, indeed they are usually quite large.

A recent realization result has been obtained by Brustenga and the
author in [\refcite{AB}]. As we will note later, the two results
mentioned above (Theorem \ref{thm:dim} and Theorem \ref{thm:purinf})
can be seen as particular cases of the main result
in [\refcite{AB}]. Before we proceed with the statement of this result,
and in order to put it in the right setting, we need to recall a few
monoid theoretic concepts.

\begin{definition}
\label{prime} {\rm Let $M$ be a monoid. An element $p\in M$ is {\it prime}
if for all $a_1,a_2\in M$, $p\le a_1+a_2$ implies $p\le a_1$ or
$p\le a_2$. A monoid is {\it primely generated} if each of its elements is
a sum of primes.}
\end{definition}

\begin{proposition} \cite[Corollary 6.8]{Brook01}
\label{fingen} Any finitely generated refinement monoid is primely
generated.
\end{proposition}

We have the following particular case of question R1.

\bigskip

\noindent {\bf R2. Realization Problem for finitely generated
refinement monoids:} Is every finitely generated conical refinement
monoid realizable by a von Neumann regular ring?

\bigskip

We conjecture that R2 has a positive answer.  Our main tool to
realize a large class of finitely generated refinement monoids is
the consideration of some regular algebras associated with quivers.

\bigskip

Recall that a  {\it quiver} (= {\it directed graph}) consists of a `vertex
set' $E^0$, an `edge set' $E^1$, together with maps~$r$ and~$s$
from $E^1$ to $E^0$ describing, respectively, the range and source
of edges. A quiver $E=(E^0,E^1,r,s)$ is said to be {\it row-finite}
in case, for each vertex $v$, the set
$s^{-1}(v)$ of arrows with source $v$ is finite.
For a row-finite quiver $E$, the graph monoid $M(E)$ of $E$ is defined as the
quotient monoid of $F=F_E$, the free abelian monoid with basis $E^0$
modulo the congruence generated by the relations
$$v=\sum _{\{ e\in E^1\mid s(e)=v \}} r(e) $$
for every vertex $v\in E^0$ which emits arrows (that is
$s^{-1}(v)\ne \emptyset$).

It follows from Proposition \ref{fingen} and \cite[Proposition 4.4]{AMP} that, for a finite quiver
$E$, the monoid $M(E)$ is primely generated. Note that this is not
always the case for a general row-finite graph $E$. An example is
provided by the graph:

$$
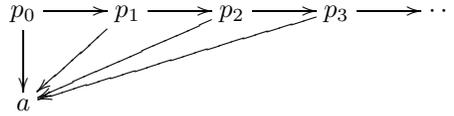
\xymatrix{p_0\ar[r] \ar[d]  & p_1 \ar[r]\ar[dl] & p_2 \ar[r]\ar[dll] & p_3 \ar[r]\ar[dlll] & \cdots \\
a}$$

The corresponding monoid $M$ has generators $a,p_0,p_1, \dots $
and relations given by $p_i=p_{i+1}+a$ for all $i\ge 0$. One can
easily see that the only prime element in $M$ is $a$, so that $M$
is not primely generated.

Now we have the following result of Brookfield:

\begin{theorem}
\label{primgensep} \cite[Theorem 4.5 and Corollary 5.11(5)]{Brook01}
Let $M$ be a primely generated refinement monoid. Then $M$ is
separative.
\end{theorem}

In fact, primely generated refinement monoids enjoy many other nice
properties, see [\refcite{Brook01}] and also [\refcite{WDim}].

It follows from Proposition \ref{fingen} and Theorem
\ref{primgensep} that a finitely generated refinement monoid is
separative. In particular all the monoids associated to finite
quivers are separative. For a row-finite quiver $E$ the result
follows by using the fact that the monoid $M(E)$ is the direct limit of
monoids associated to certain finite subgraphs of $E$, see
\cite[Lemma 2.4]{AMP}.

\medskip

\begin{theorem}
\label{thm:quivers} {\rm (\cite[Theorem 4.2, Theorem 4.4]{AB})}
Let $M(E)$ be the monoid corresponding to a finite quiver $E$ and let $K$
be any field. Then there exists a unital von Neumann
regular hereditary $K$-algebra $Q_K(E)$ such that $\mon{Q_K(E)}\cong M(E)$.
Furthermore, if $E$ is a row-finite quiver, then there exists a
(not necessarily unital) von Neumann regular $K$-algebra $Q_K(E)$ such that
$\mon{Q_K(E)}\cong M_E$.
\end{theorem}

\medskip

Note that, due to
unfortunate lack of convention in this area, the Cuntz-Krieger relations used in
[\refcite{AB}] are the opposite to the ones used in [\refcite{AMP}], which are the ones we
are following in this survey, so that the result in \cite[Theorem 4.4]{AB}
is stated for {\it column-finite} quivers instead of row-finite ones.
The regular algebras $Q_K(E)$ are related to the {\it Leavitt path algebras}
$L_K(E)$ of [\refcite{AA1}], [\refcite{AA2}], [\refcite{AMP}].

We now observe that Theorem \ref{thm:dim} and Theorem \ref{thm:purinf}
are particular cases of Theorem \ref{thm:quivers}.
This follows from the fact that the monoids considered in these theorems
are known to be graph monoids $M(E)$ for suitable quivers $E$.
Indeed, taking into account \cite[Theorem 3.5 and Theorem 7.1]{AMP}, we see that the case of
dimension monoids follows from \cite[Proposition 2.12]{Raeburn} and
the case of monoids of the form $\{ 0 \}\sqcup G$, with $G$ a
countable abelian group, follows from \cite[Theorem 1.2]{Szy}.
As Pardo pointed out to me, the quiver $E$ can be chosen to be finite in case
$M=\{0\}\sqcup G$ for a finitely generated abelian group $G$. To see this, note that
such a monoid admits a presentation given by a finite number of generators $a_1,\dots ,a_n$,
and relations of the form $a_i=\sum _{j=1}^n \gamma_{ji}a_j$, where all $\gamma _{ji}$
are strictly positive integers and $\gamma _{ii}\ge 2$ for all $i$.
The corresponding finite quiver will have $\gamma _{ji}$ arrows from the vertex $i$
to the vertex $j$.

Now we would like to describe how  this construction sheds light on
problem R2. The answer is completely known for the class of
antisymmetric finitely generated refinement monoids. The monoid $M$
is said to be {\it antisymmetric} in case the algebraic pre-order is
a partial order, that is, in case $x\le y$ and $y\le x$ imply that
$x=y$. Note that every antisymmetric monoid is conical.

We say that a monoid $M$ is {\it primitive} if it is an
antisymmetric primely generated refinement monoid \cite[Section
3.4]{Pierce}. A primitive monoid $M$ is completely determined by its set
of primes $\mathbb P (M)$ together with a transitive and
antisymmetric relation $\lhd$ on it, given by $q\lhd p$ iff $p+q=p$.
Indeed given such a pair $(\mathbb P, \lhd )$, the abelian monoid
$M(\mathbb P,\lhd)$ defined by taking as a set of generators
$\mathbb P$ and with relations given by $p=p+q$ whenever $q\lhd p$,
is a primitive monoid, and the correspondences $M\mapsto (\mathbb P
(M),\lhd) $ and $(\mathbb P,\lhd)\mapsto M(\mathbb P,\lhd)$ give a
bijection between isomorphism types of primitive monoids and
isomorphism types of pairs $(\mathbb P,\lhd)$, where $\mathbb P$ is
a set and $\lhd $ a transitive antisymmetric relation on $\mathbb
P$, see \cite[Proposition 3.5.2]{Pierce}.

Let $M$ be a primitive monoid and $p\in \mathbb P (M)$. Then $p$ is
said to be {\it free} in case $p\ntriangleleft p$. Otherwise $p$ is
{\it regular}, see \cite[Section 2]{APW}. So giving a primitive
monoid is equivalent to giving a poset $(\mathbb P,\le )$ which is a
disjoint union of two subsets: $\mathbb P=\Pfree \sqcup \Preg $. If
$M$ is a finitely generated primitive monoid then $\mathbb P (M)$ is
a finite set, indeed $\mathbb P (M)$ is the minimal generating set
of $M$.


We can now describe the finitely generated primitive monoids which are graph monoids.
Recall that a {\it lower cover} of an element $p$ of a poset $\mathbb P$
is an element $q$ in $\mathbb P$ such that $q<p$ and
$[q,p]=\{q,p\}$. The set of lower covers of $p$ in $\mathbb P (M)$ is denoted by
$\rL(M,p)$, and $\rLfree(M,p)$ and $\rLreg(M,p)$ denote the sets of
free and regular elements in $\rL(M,p)$ respectively.

\begin{theorem}\label{T:charact} {\rm (cf. \cite[Theorem 5.1]{APW})}
Let $M$ be a finitely generated primitive monoid. Then the following
statements are equivalent:
\begin{enumerate}
\item $M$ is a graph monoid.

\item $M$ is a direct limit of graph monoids.

\item $|\rLfree(M,p)|\leq 1$ for each $p\in \Pfree(M)$.
\end{enumerate}
\end{theorem}

For the monoids as in the statement there is a hereditary, von Neumann
regular ring $Q(E)$
such that $\mon{Q(E)}=M_E=M$ (Theorem \ref{thm:quivers}).
It is worth mentioning that, in some
cases, an infinite quiver is required in Theorem \ref{T:charact}.
A slightly more general result is indeed presented in \cite[Theorem
5.1]{APW}. Namely the same characterization holds when $M$ is a
primitive monoid such that $\rL(M,p)$ is finite for every $p$ in
$\mathbb P (M)$.

In view of Theorem \ref{T:charact}, the simplest primitive monoid which
is not a graph monoid is the monoid
$$M=\langle p,a,b\mid p=p+a=p+b \rangle .$$

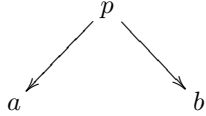
\begin{figure}[htb]
 \[
 {
 \xymatrix{
 & p \ar[ld]\ar[rd] & \\
 a & & b
 }
 }
 \]
\caption{The poset $\mathbb P (M)$ for the monoid $M=\langle
p,a,b\mid p=p+a=p+b \rangle $.} \label{Fig:Graphpabp}
\end{figure}

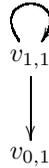
\begin{figure}[htb]
 \[
 {
 \xymatrix{
 & v_{1,1}\ar@(ul,ur)\ar[d]  \\
 & v_{0,1}
 }
 }
 \]
\caption{The graph $S_1$.} \label{Fig:QuiverS1}
\end{figure}

In this example $\mathbb P (M)=\Pfree (M)=\{p,a,b\}$, and $p$ has
two free lower covers $a,b$. Thus, by Theorem \ref{T:charact}, the monoid $M$
is not even a direct limit of graph monoids (with monoid homomorphisms as
connecting maps). However $M$ can be realized as the monoid of a
suitable von Neumann regular ring, as follows. Fix a field $K$ and consider two
indeterminates $t_1,t_2$ over $K$. We consider the regular algebra
$Q_{K(t_2)}(S_1)$ over the quiver $S_1$ with two vertices $v_{0,1},
v_{1,1}$ and two arrows $e_1,f_1$ such that
$r(e_1)=s(e_1)=v_{1,1}=s(f_1)$ and $r(f_1)=v_{0,1}$. The picture of
$S_1$ is as shown in Figure \ref{Fig:QuiverS1}.

The algebra $Q_1:=Q_{K(t_2)}(S_1)$ has a unique non-trivial
(two-sided) ideal $M_1$, which coincides with its socle, so that we
get an extension of rings:
$$
\begin{CD}
0 @>>> M_1 @>>> Q_1   @>\pi_1>>  K(t_2)(t_1)=K(t_1,t_2) @>>> 0
\end{CD}
$$
However the element $t_1$ does not lift to a unit in $Q_1$, rather
there are $z_1,\ol{z}_1$ in $Q_1$ such that $\ol{z}_1z_1=1$, but
$z_1\ol{z}_1\ne 1$, and $\pi _1(z_1)=t_1$. Some additional information on the
algebra $Q_F(S_1)$, where $F$ denotes an arbitrary field, can be found in \cite[Examples 4.3]{AB}.

Let $S_2$ be a copy of $S_1$, now with vertices labelled as
$v_{0,2}, v_{1,2}$ and arrows labelled as $e_2,f_2$, and set
$Q_2=Q_{K(t_1)}(S_2)$. There is a corresponding diagram
$$
\begin{CD}
0 @>>> M_2 @>>> Q_2   @>\pi_2>>  K(t_1)(t_2)=K(t_1,t_2) @>>> 0
\end{CD}
$$

Let $P$ be the pullback  of the maps $\pi_1$ and $\pi _2$, so that
$P$ fits in the following commutative square:
\begin{equation}
\label{equ:4.1}
\begin{CD}
P @>\rho _1>> Q_1\\
   @V\rho_2 VV  @VV\pi _1V\\
Q_2 @>\pi_2>> K(t_1,t_2)
\end{CD}
\end{equation}
Then $P$ is a von Neumann regular ring and $\mon{P}=M$, see [\refcite{Aposet}]. Indeed a
wide generalization of this method gives the following realization result:

\begin{theorem}
\label{primitive} {\rm (\cite[Theorem 2.2]{Aposet})}
Let $M$ be a finitely generated primitive monoid such that all
primes of $M$ are free and let $K$ be a field. Then there is a unital regular $K$-algebra
$Q_K(M)$ such that $\mon{Q_K(M)}\cong M$.
\end{theorem}

Moreover both the regular algebras associated with quivers [\refcite{AB}]
and the regular algebras constructed in [\refcite{Aposet}] are given explicitly in terms
of generators and relations (including universal localization [\refcite{scho}]).

Recall that a monoid $M$ is
{\it strongly separative} in case $a+a=a+b$ implies $a=b$ for
$a,b\in M$. A ring $R$ is said to be {\it strongly separative}
in case $\mon{R}$ is a strongly separative monoid, see [\refcite{AGOP}] for
background and various equivalent conditions. As we mentioned above,
every primely generated refinement monoid is separative \cite[Theorem
4.5]{Brook01}. In particular every primitive monoid is separative.
Moreover, a primitive monoid $M$ is strongly separative if and
only if all the primes in $M$ are free, see \cite[Theorem 4.5,
Corollary 5.9]{Brook01}. Thus, the class of monoids covered by
Theorem \ref{primitive} coincides exactly with the strongly separative finitely
generated primitive monoids. The case of a general finitely generated
primitive monoid remains open, although it seems amenable to analysis in the
light of [\refcite{AB}] and [\refcite{Aposet}].

\section{Realizing distributive lattices}
\label{sect:CLP}

Let $R$ be a regular ring. Then the lattice $\Id (R)$ of all
(two-sided) ideals of $R$ is an algebraic distributive lattice. Here
an {\it algebraic lattice} means a complete lattice such that each element
is the supremum of all the compact elements below it. The set of
compact elements in $\Id (R)$ is the set $\Id _c(R)$ of all finitely
generated ideals of $R$, and it is a distributive semilattice, see
for example [\refcite{GW}]. Here a {\it semilattice} means a $\vee$-semilattice with least
element $0$. A semilattice is the same as a monoid $M$
such that $x=x+x$ for every $x$ in $M$, and a distributive semilattice is
just a semilattice satisfying the refinement axiom \cite[Lemma 2.3]{GW}.
Observe that, if $M$ is a semilattice then
$x\vee y=x+y$ gives the supremum of $x,y$ in $M$.

\vspace{.2cm}

The famous {\it Congruence Lattice Problem} (CLP) asks whether an
algebraic distributive lattice is the congruence lattice of some
lattice; equivalently, whether every distributive semilattice is the
semilattice of all the compact congruences of a lattice. This
problem has been recently solved in the negative by Fred Wehrung
[\refcite{WAdv}], who constructed for each $\aleph \ge \aleph_{\omega
+1}$ an algebraic distributive lattice with $\aleph$ compact
elements such that it cannot be represented as the congruence
lattice of any lattice, see also [\refcite{GratNotices}].
His methods have been refined by R\r{u}\v{z}i\v{c}ka [\refcite{Ruz}] to
cover the case $\aleph \ge \aleph _2$. It is worth to remark that
Wehrung had previously shown in [\refcite{WPub}] that every algebraic distributive
lattice with $\le \aleph _1$ compact elements can be realized as
the ideal lattice for some von Neumann regular ring $R$, and thus
is isomorphic to the congruence lattice of the lattice $ L
(R_R)$ of principal right ideals of $R$ \cite[Corollary 4.4]{WPams}.
So the formulation of R\r{u}\v{z}i\v{c}ka is best possible
(concerning the number of compact elements).
The ``hard core" of Wehrung's proof in [\refcite{WPub}] is a ring-theoretic amalgamation
result proved by P. M. Cohn in \cite[Theorem 4.7]{Cohn}.

\vspace{.2cm}

One can ask: What is the relationship between the CLP, or more
concretely, the representation problem of algebraic distributive
lattices as lattices of ideals of regular rings, and our problem
R1? The answer is that, for a regular ring $R$, the
lattice $\Id (R)$ is only a small piece of information compared with
the information contained in the monoid $\mon{R}$, in the sense that
$\mon{R}$ determines $\Id (R)$, but generally the structure of
$\mon{R}$ can be much more complicated than the structure of $\Id
(R)$, e.g. for simple rings. Indeed we have a lattice isomorphism
$\Id (R)\cong \Id (\mon{R})$, where for a conical monoid $M$, $\Id
(M)$ is the lattice of all order-ideals of $M$, cf.
\cite[Proposition 7.3]{GW}.  Recall that an
{\it order-ideal} of $M$ is a submonoid $I$ of $M$ with the property that
whenever $x\le y$ in $M$ and $y\in I$, we have $x\in I$.
If $L$ is an algebraic distributive lattice which is not the congruence
lattice of any lattice and $M$ is any conical refinement monoid such that
$\Id (M)\cong L$, then $M$ cannot be realized as $\mon{R}$ for a regular ring $R$.
For every algebraic distributive lattice $L$ there is at least one such conical refinement
monoid, namely the semilattice $L_c$ of compact elements of $L$, but we should
expect a myriad of such monoids to exist. Wehrung proved in [\refcite{WPub}]
that if $|L_c|\le\aleph _1$ then the semilattice $L_c$ can be realized
as $\mon{R}$ where $R$ is a von Neumann regular ring, and he showed in [\refcite{Wtams}]
that there is a distributive semilattice $S_{\omega_1}$  of size
$\aleph _1$ which is not the semilattice of finitely generated, idempotent-generated
ideals of any exchange ring of finite stable rank. In particular there is no
locally matricial $K$-algebra $A$ over a field $K$ such that $\Id_c(A)\cong S_{\omega _1}$; see
[\refcite{Wtams}] for details. This contrasts with Bergman's result [\refcite{Bergman}]
stating that every distributive semilattice of size $\le \aleph _0$ is the semilattice of finitely generated
ideals of an ultramatricial $K$-algebra, for every field $K$.

Say that a subset $A$ of a poset $\mathbb P$ is a {\it lower
subset} in case $q\le p$ and $p\in A$ imply $q\in A$. The
set $\mathcal L (\mathbb P)$ of all lower subsets of $\mathbb P$
forms an algebraic distributive lattice, which is a sublattice of
the Boolean lattice ${\bf 2}^{\mathbb P}$.
Now if $L$ is a {\it finite} distributive lattice, then by a result
of Birkhoff (\cite[Theorem II.1.9]{Grat}) there is a finite poset $\mathbb P$ such that $L$ is the
lattice of all lower subsets of $\mathbb P$. In the case of a finite
Boolean algebra ${\bf 2}^n$ with $n$ atoms, the poset $\mathbb P$ is just an
antichain with $n$ points, and ${\bf 2}^n=\mathcal P (\mathbb P)=\mathcal L(\mathbb P)$. Our
construction in [\refcite{Aposet}] gives a realization of $L$
as the ideal lattice of a regular $K$-algebra  $Q_K(\mathbb P)$,
where $K$ is an arbitrary fixed field,
such that the monoid $\mon{Q_K(\mathbb P)}$ is the monoid $M(\mathbb P)$
associated with $(\mathbb P,<)$, with all elements
in $M(\mathbb P)$ being free, that is, $M(\mathbb P)$ is the abelian monoid
with generators $\mathbb P$
and relations given by $p=q+p$ whenever $q<p$ in $\mathbb P$.
Moreover $Q_K(\mathbb P)$ satisfies the following properties
(\cite[Proposition 2.12, Remark 2.13, Theorem 2.2]{Aposet}):
\begin{enumerate}[(a)]
\item
There is a canonical family of commuting idempotents $\{e(A):A\in
L \}$ such that
\begin{enumerate} [(i)]
\item $e(A)e(B)=e(A\cap B)$
\item $e(A)+e(B)-e(A)e(B)=e(A\cup B)$
\item $e(\emptyset)=0$ and $e(\mathbb P)=1$.
\item $e(A)Q_K(\mathbb P)e(A)\cong Q_K(A)$.
\end{enumerate}
\item Let $I(A)$ be the ideal of $Q_K(\mathbb P)$ generated by $e(A)$.
Then the assignment $$A\mapsto I(A)$$
defines a lattice isomorphism from $L=\mathcal L(\mathbb P)$ onto
$\Id (Q_K(M))$.
\item The map $M(\mathbb P)\to \mon{Q_K(\mathbb P)}$ given by $p\mapsto [e(\mathbb P\dnw p)]$,
for $p\in \mathbb P$,
is a monoid isomorphism. Here $\mathbb P\dnw p=\{q\in \mathbb  P:q\le p\}$ is the lower subset of $\mathbb P$
generated by $p$.
\end{enumerate}

\noindent The set $\text{Idem}(R)$ of idempotents of a ring is a poset in a natural way, by using
the order $e\le f$ iff $e=fe=ef$. This poset is a partial lattice, in the following sense: every
two commuting idempotents $e$ and $f$ have an infimum $ef$ and a supremum $e+f-ef$ in $\text{Idem}(R)$.
 Say that a map $\phi \colon L\to \text{Idem}(R)$ from a lattice $L$
to $\text{Idem}(R)$ is a {\it lattice homomorphism} in case $\phi (x)$ and $\phi (y)$
commute and $\phi (x\vee y)=\phi (x)\vee \phi (y)$ and $\phi (x\wedge y)=\phi (x)\wedge \phi (y)$,
for every $x,y\in L$.

The above results can be paraphrased as follows:
The canonical mapping $\text{Idem}(Q_K(\mathbb P))\to \Id(Q_K(\mathbb P))=\mathcal L(\mathbb P)$
sending each element
$e$ in $\text{Idem}(Q_K(\mathbb P))$ to the ideal generated by $e$ has a distinguished
section $e\colon \Id(Q_K(\mathbb P))\to \text{Idem}(Q_K(\mathbb P))$, $A\mapsto e(A)$, which is
a lattice homomorphism.

Write $Q=Q_K(\mathbb P)$. Observe that we have lattice isomorphisms
$$L=\mathcal L(\mathbb P)\cong\Id (Q)\cong \Id \mon{Q}\cong \mon{Q}/{\asymp} .$$

Here $\mon{Q}/{\asymp}$ is the maximal semilattice quotient of the monoid $\mon{Q}$,
see \cite[Section 2]{GW}.
The finite distributive lattice $L$ can be represented in many other
ways as an ideal lattice of a regular ring, for instance using
ultramatricial algebras [\refcite{Bergman}], but the monoids corresponding to these
ultramatricial algebras have little to do with $M(\mathbb P)$. Indeed as soon as $\mathbb P$
is not an antichain we will have that $\mon{R}$ is non-finitely generated for every ultramatricial
algebra $R$ such that $\Id (R)\cong \mathcal L (\mathbb P)$.

\section{The dependence on the field}
\label{sect:field}

We like to work with von Neumann regular rings which are algebras
over a field $K$. A natural question is whether the field $K$ plays
any role concerning the realization problem. So we ask the following
variant of R1.

\medskip

\noindent {\bf $\text{R}(K)$. Realization Problem for von Neumann Regular $K$-algebras}
Let $K$ be a fixed field. Is every countable conical refinement
monoid realizable by a von Neumann regular $K$-algebra?

\medskip

The answer to this question is known for uncountable fields, thanks
to an observation due to Wehrung. Indeed the basic counter-example
comes from a construction due to Chuang and Lee [\refcite{CL}].
Their remarkable example gave a negative answer to five open
questions in Goodearl's book [\refcite{vnrr}] on von Neumann regular
rings.

We take this opportunity to present the complete argument, including
a result of Goodearl, generalizing  Wehrung's observation, and a
version of the Chuang-Lee construction.

We will need the notion of (pseudo-)rank function, as given in
\cite[Chapter 16]{vnrr}. Recall that a {\it pseudo-rank function}
$N$ on a unital regular ring $R$ is a function $N\colon R\to [0,1]$
such that
\begin{enumerate} [(a)]
\item $N(1)=1$.
\item $N(xy)\le N(x)$ and $N(xy)\le N(y)$ for all $x,y\in R$.
\item $N(e+f)=N(e)+N(f)$ for all orthogonal idempotents $e,f\in R$.
\end{enumerate}
A {\it rank function} is a pseudo-rank function $N$ such that
$N(x)>0$ for all nonzero $x\in R$.

\begin{proposition}\label{norepF} (Goodearl)
Let $(M,u)$ be a conical refinement monoid with order-unit admitting
a faithful state $s$, i.e. a monoid homomorphism $s\colon M\to \R^+$
such that $s(u)=1$ and $s(x)>0$ for every nonzero $x$ in $M$. Assume
that $M$ is not cancelative. Then there is no regular algebra $R$
over an uncountable field $F$ such that $\mon{R}\cong M$.
\end{proposition}

\begin{proof}
Assume that $R$ is a regular $F$-algebra over an uncountable field
$F$ with $\mon{R}\cong M$. Clearly we can assume that $R$ is unital
and that $[1]$ corresponds to $u$ under the isomorphism
$\mon{R}\cong M$.

By \cite[Theorem 2.2]{GM}, it suffices to prove that there is no
uncountable independent family of nonzero right or left
ideals of $R$. Since $\mon{R}\cong \mon{R^{\text{op}}}$, where $R^{\text{op}}$
is the opposite ring of $R$, we see that it suffices to show this
fact for right ideals. Indeed, once this is established, we get from
\cite[Theorem 2.2]{GM} and \cite[Proposition 4.12]{AGOP} that $R$ is unit-regular,
and thus $\mon{R}$ must be cancellative by \cite[Theorem 4.5]{AGOP},
a contradiction with our hypothesis.

By \cite[Proposition 17.12]{vnrr} there exists a pseudo-rank
function $N$ on $R$ such that $N(x)=s([xR])$ for every $x\in R$.
Since $s$ is faithful, we see that $N$ is indeed a rank function.
Now by \cite[Proposition 16.11]{vnrr} we get that $R$ contains no
uncountable direct sums of nonzero right ideals, as desired.
\end{proof}

Now we are going to recall the example of Chuang and Lee [\refcite{CL}].
We will give a presentation which is a little bit more general. Let
$R$ be a $\sigma$-unital regular ring, that is, a regular ring
having an increasing sequence $(e_n)$ of idempotents in $R$ such
that $R=\bigcup _{n=1}^{\infty}e_nRe_n$. Put $R_n=e_nRe_n$. Recall that
the multiplier ring $\mathcal M (R)$ is the completion of $R$ with
respect to the strict topology; see [\refcite{APCom}]. Write
$$\mathcal R=\{(x_n)\mid (x_n) \text{ is a Cauchy sequence in the strict topology } \}
\subseteq \prod_{n=1}^{\infty} R_n .$$  By the continuity of
operations, $\mathcal R$ is a unital subring of
$\prod_{n=1}^{\infty} R_n $. There is an obvious canonical
surjective homomorphism $\Phi \colon \mathcal R\to \mathcal M (R)$
whose kernel is $I=\{(x_n)\mid x_n\to 0\}$, where the convergence is
with respect to the strict topology.

\begin{lemma}
\label{Ireg} $I$ is always a (non-unital) regular ring. If each $e_nRe_n$ is
unit-regular, then $I$ is unit-regular, meaning that $eIe$ is
unit-regular for every idempotent $e$ in $I$.
\end{lemma}

\begin{proof}
Let $x=(x_n)\in I$. Choose a sequence of integers $m_1<m_2<\cdots $
such that for all $m\ge m_i$ we have $x_me_i=e_ix_m=0$. Now for
$m_i\le m<m_{i+1}$, choose a quasi-inverse $y_m$ of $x_m$ in
$(e_m-e_i)R(e_m-e_i)$. (Note that $x_m\in (e_m-e_i)R(e_m-e_i)$ for
$m_i\le m<m_{i+1}$.) We get a quasi-inverse $y=(y_n)$ of $x$ such
that $y_n\to 0$ strictly, so $y\in I$ and $I$ is regular. The last
part is easy, and is left to the reader.
\end{proof}

Observe that if $Q$ is any regular ring such that $Q\subseteq
\mathcal M (R)$, then $\Phi ^{-1}(Q)$ is a regular ring (\cite[Lemma
1.3]{vnrr}) which is a subdirect product of the regular rings
$(R_n)$. In particular $\Phi ^{-1}(Q)$ is stably finite if each
$R_n$ is so.

Now we see that when $K$ is a {\it countable field} the regular
algebra $Q_K(E)$ of the quiver $E$ with $E^0=\{v_0,v_1\}$ and
$E^1=\{e,f\}$, with $r(e)=s(e)=s(f)=v_1$ and $r(f)=v_0$ gives an
example that fits in the above picture. (Note that the quiver $E$ is
the same as the quiver $S_1$ of Figure  \ref{Fig:QuiverS1}.)
Since the field $K$ is countable, the algebra $Q_K(E)$ is also countable.

Write $Q=Q_K(E)$, and let $I$ be the ideal of $Q$ generated by $v_0$. Then
$I=\text{Soc}(Q)$ is a simple (non-unital) ring, and $I$ is countable,
so we have $I\cong M_{\infty}(K)$ (because $v_0Qv_0$ is isomorphic to $K$),
see \cite[Remark 2.9]{APCom}.
Here $M_{\infty}(K)$ denotes the $K$-algebra of countably infinite matrices with only a
finite number of nonzero entries. Since this is a crucial argument here,
let us recall the details. The ring $I$ is countable and simple with a minimal
idempotent $v_0$, so by general theory there is a dual pair $V, W$
of $K$-vector spaces such that $I\cong \mathcal F _W(V)$, the algebra of all
adjointable operators on $V$ of finite rank. Since $I$ is countable, both $V$
and $W$ are countably dimensional $K$-vector spaces. By an old result of G. W. Mackey
\cite[Lemma 2]{Mackey}, there are dual bases $(v_i)$ and $(w_j)$
for $V$ and $W$ respectively, that is, we have
$\langle v_i,w_j\rangle =\delta _{ij}$ for all $i,j$, which shows that
$\mathcal F_W(V)\cong M_{\infty}(K)$. Since $I$ is essential in $Q$ we
get an embedding of $Q$ into the multiplier algebra $\mathcal M
(I)$. Observe that $\mathcal M (I)\cong RCFM(K)$, the algebra of
row-and-column-finite matrices with coefficients in $K$, and that
$Q/I\cong K(t)$. So the algebra $Q$ has the same essential
properties as the Chuang and Lee algebra, see [\refcite{CL}]. Now
$M_{\infty}(K)$ is clearly $\sigma$-unital and unit-regular. Indeed
there is a $\sigma$-unit $(e_n)$ for $I$ consisting of idempotents
such that $e_nIe_n\cong M_n(K)$. Now Lemma \ref{Ireg} together with
\cite[Lemma 1.3]{vnrr} gives that $S:=\Phi^{-1}(Q)$ is regular, and
it is residually artinian. The ring $S$ is not countable but it can
be easily modified to get a countable algebra with similar
properties. Indeed consider the $K$-subalgebra $S_0$ of $S$
generated by $\oplus _{n=1}^{\infty}M_n(K)$ and $a,b$ where $a,b$
are elements in $S$ such that $\Phi (a)\Phi (b)=1$ and $\Phi (b)\Phi
(a)\ne 1$. Observe that $S_0$ is countable. We can build a sequence
of countable $K$-subalgebras of $S$:
$$S_0\subseteq S_1\subseteq S_2\subseteq \cdots \subseteq S$$
such that each element in $S_i$ is regular in $S_{i+1}$ for all $i$.
It follows that $S_{\infty}=\bigcup S_i$ is a countable, regular
$K$-algebra, which is embedded in $\prod _{n=1}^{\infty}M_n(K)$.
Moreover $S_{\infty}$ cannot be unit-regular because it has a
quotient ring which is not directly finite. It follows that
$M=\mon{S_{\infty}}$ is not cancellative and it is a countable
monoid satisfying the hypothesis of Proposition \ref{norepF},
because there is a rank function on $\prod _{n=1}^{\infty}M_n(K)$,
e.g. the function $N((x_n))=\sum _{n=1}^{\infty} 2^{-n}N_n(x_n)$,
where $N_n$ is the unique rank function on $M_n(K)$. Therefore $M$
gives a counterexample to $\text{R}(F)$ for uncountable fields $F$,
although by definition it can be realized over {\it some} countable
field $K$.

\section*{Acknowledgments}

This work has been partially supported by the DGI
and European Regional Development Fund, jointly, through Project
MTM2005-00934, and by the Comissionat per Universitats i Recerca
de la Generalitat de Catalunya.

It is a pleasure to thank Gene Abrams, Ken Goodearl, Kevin O'Meara and
Enrique Pardo for their helpful comments. I am specially grateful to
Fred Wehrung for his many valuable comments and suggestions.

\end{document}